\documentclass[a4paper, 11pt]{amsart}

\usepackage[colorlinks=true,citecolor=green,linkcolor=blue]{hyperref}
\usepackage{amsmath}
\usepackage{verbatim}
\usepackage{amsfonts}
\usepackage{amssymb}
\usepackage{blkarray}
\usepackage{enumerate}
\usepackage[inner=2.8cm,outer=2.8cm, bottom=3.5cm]{geometry}
\usepackage{mathrsfs}
\usepackage{stmaryrd}
\usepackage{setspace}
\usepackage{bbold}
\renewcommand{\k}{\mathbb{k}}

\usepackage{tikz}
\usetikzlibrary{matrix}

\usepackage{aliascnt}
\usepackage[initials]{amsrefs}

\newcommand{\bbM}{\mathbb{M}}

\newcommand{\NN}{\normalfont\mathbb{N}}
\newcommand{\ZZ}{\normalfont\mathbb{Z}}

\newcommand{\PP}{\normalfont\mathbb{P}}

\newcommand{\mm}{{\normalfont\mathfrak{m}}}

\newcommand{\pp}{{\normalfont\mathfrak{p}}}

\newcommand{\rank}{\normalfont\text{rank}}

\newcommand{\depth}{\normalfont\text{depth}}
\newcommand{\grade}{\normalfont\text{grade}}

\newcommand{\HT}{\normalfont\text{ht}}

\newcommand{\Sym}{\normalfont\text{Sym}}
\newcommand{\Rees}{\mathcal{R}}

\newcommand{\Fitt}{\normalfont\text{Fitt}}
\newcommand{\EEQ}{\mathcal{K}}

\newcommand{\OO}{\mathcal{O}}

\newcommand{\LL}{\mathbb{L}}

\newcommand{\FF}{\normalfont\mathfrak{F}}
\newcommand{\sF}{\widetilde{\mathfrak{F}}}
\newcommand{\HL}{\normalfont\text{H}_{\mm}}
\newcommand{\HH}{\normalfont\text{H}}

\newcommand{\AAA}{\mathcal{A}}

\newcommand{\bideg}{\normalfont\text{bideg}}
\newcommand{\Proj}{\normalfont\text{Proj}}

\newtheorem{headthm}{Theorem}

\newaliascnt{headcor}{headthm}
\newtheorem{headcor}[headcor]{Corollary}
\aliascntresetthe{headcor}

\newaliascnt{headconj}{headthm}

\aliascntresetthe{headconj}

\newaliascnt{corollary}{theorem}

\aliascntresetthe{corollary}

\newaliascnt{lemma}{theorem}
\newtheorem{lemma}[lemma]{Lemma}
\aliascntresetthe{lemma}

\newaliascnt{conjecture}{theorem}

\aliascntresetthe{conjecture}

\newaliascnt{proposition}{theorem}
\newtheorem{proposition}[proposition]{Proposition}
\aliascntresetthe{proposition}

\theoremstyle{definition}
\newaliascnt{definition}{theorem}
\newtheorem{definition}[definition]{Definition}
\aliascntresetthe{definition}

\newaliascnt{notation}{theorem}

\aliascntresetthe{notation}

\newaliascnt{example}{theorem}

\aliascntresetthe{example}

\newaliascnt{examples}{theorem}

\aliascntresetthe{examples}

\newaliascnt{remark}{theorem}
\newtheorem{remark}[remark]{Remark}
\aliascntresetthe{remark}

\newaliascnt{problem}{theorem}

\aliascntresetthe{problem}

\newaliascnt{construction}{theorem}

\aliascntresetthe{construction}

\newaliascnt{setup}{theorem}
\newtheorem{setup}[setup]{Setup}
\aliascntresetthe{setup}

\newaliascnt{algorithm}{theorem}

\aliascntresetthe{algorithm}

\newaliascnt{observation}{theorem}

\aliascntresetthe{observation}

\newaliascnt{defprop}{theorem}

\aliascntresetthe{defprop}

\def\equationautorefname~#1\null{(#1)\null}
\def\sectionautorefname~#1\null{Section #1\null}
\def\subsectionautorefname~#1\null{\S #1\null}

\newcommand{\Spec}{\operatorname{Spec}}

\newcommand{\D}{\mathcal{D}}

\renewcommand{\leq}{\leqslant}

\title[Saturated special fiber ring of height three Gorenstein ideals]{Multiplicity of the saturated special fiber ring of height three Gorenstein ideals}

\date{\today}

\author{Yairon Cid-Ruiz}
\address[Cid-Ruiz]{Max Planck Institute for Mathematics in the Sciences, Inselstra\ss e 22, 04103 Leipzig, Germany.}
\email{cidruiz@mis.mpg.de}
\urladdr{https://ycid.github.io}

\author{Vivek Mukundan}
\address[Mukundan]{Department of Mathematics, University of Virginia, Charlottesville, VA 22904-4135, USA.}
\email{vm6y@virginia.edu}


\begin{document}

\begin{abstract}
	Let $R$ be a polynomial ring over a field and $I \subset R$ be a Gorenstein ideal of height three that is minimally generated by homogeneous polynomials of the same degree.
	We compute the multiplicity of the \textit{saturated special fiber ring} of $I$.
	The obtained formula depends only on the number of variables of $R$, the minimal number of generators of $I$, and the degree of the syzygies of $I$.
	Applying results from \cite{MULTPROJ}, we get a formula for the $j$-multiplicity of $I$ and an effective method to study a rational map determined by a minimal set of generators of $I$.
\end{abstract}

\maketitle

\section{Introduction}

The \textit{saturated special fiber ring} is an algebra that was in introduced in \cite{MULTPROJ} and that was initially motivated by the interest of studying rational maps with the use of syzygies and blow-up algebras.   
This algebra encodes valuable information regarding the degree and birationality of rational maps and is also related to the $j$-multiplicity of ideals (see \cite{MULTPROJ}).
Recently, for the case of perfect ideals of height two, its multiplicity was computed in \cite{MULT_SAT_PERF_HT_2}. 
In this paper, we extend the family of ideals for which the multiplicity of the saturated special fiber ring is known. 
Specifically, we treat the case of Gorenstein ideals of height three that are minimally generated by homogeneous polynomials of the same degree.

\medskip

Let $\k$ be a field, $R$ be the polynomial ring $R=\k[x_1,\ldots,x_d]$ and $\mm \subset R$ be the maximal irrelevant ideal $\mm=(x_1,\ldots,x_d)$.
Let $I \subset R$ be a homogeneous Gorenstein ideal of height three, $n=\mu(I)$ be the minimal number of generators of $I$, and suppose that $I$ is minimally generated by homogeneous polynomials of the same degree $\delta$.
From the Buchsbaum-Eisenbud structure theorem \cite{BuchsbaumEisenbud}, we can assume that $I$ has an alternating minimal presentation matrix $\varphi$ whose non-zero entries are all of the same degree $D \ge 1$ and that $I$ is minimally generated by the $(n-1)\times (n-1)$ pfaffians of $\varphi$.
Accordingly, we can assume that $I$ is minimally generated by $n$ homogeneous polynomials of degree $\delta=\frac{1}{2}(n-1)D$.

\medskip

Following \cite{MULTPROJ}, the saturated special fiber ring of $I$ is given by the graded $\k$-algebra
$$
\sF(I) = \bigoplus_{q=0}^\infty {\left[\big(I^q:\mm^\infty\big)\right]}_{q\delta}.
$$
We also assume that the ideal $I$ satisfies the condition $G_d$, that is, $\mu(I_{\pp}) \le \dim(R_{\pp})$  for all $\pp \in V(I) \subset \Spec(R)$  such that $\HT(\pp)<d$. 
It should be noted that the condition $G_d$ is always satisfied by generic perfect ideals of height two and by generic Gorenstein ideals of height three.
We follow a strategy very similar to the one of \cite{MULT_SAT_PERF_HT_2}, which accounts to approximate the Rees algebra with the symmetric algebra.
The assumption of the condition $G_d$ allows us to the make possible that approximation.
Important tools used in this paper are: the family of complexes $\mathcal{D}_\bullet^q(\varphi)$ introduced in \cite{KU92}, and the computation \cite[Theorem 4.4]{KPU18} of the multiplicity of the special fiber ring under the further assumption that $I$ is linearly presented (i.e., when $D=1$).

\medskip

The following theorem contains the main result of this paper, where we compute the multiplicity of $\sF(I)$.

\begin{headthm}
	\label{thmA}
	Let $\k$ be a field, $R$ be the polynomial ring $R=\k[x_1,\dots,x_d]$ and $I \subset R$ be a homogeneous ideal in $R$.
	Let $n=\mu(I)$ be the minimal number of generators of $I$.
	Suppose that the following conditions are satisfied:
	\begin{enumerate}[(i)]
		\item $I$ is a height $3$ Gorenstein ideal.
		\item $n \ge d$.
		\item Every non-zero entry of an alternating minimal presentation matrix of $I$  has degree $D\ge 1$.
		\item $I$ satisfies the condition $G_d$.		
	\end{enumerate}
	Then, the multiplicity of the saturated special fiber ring $\sF(I)$ of $I$ is given by 
	$$
	e\left(\sF(I)\right) = D^{d-1} \sum_{i=0}^{\lfloor\frac{n-d}{2}\rfloor}\binom{n-2-2i}{d-2}.
	$$
\end{headthm}

Next, there are some consequences of \autoref{thmA}.

\medskip

The $j$-multiplicity of $I$ is given by
$$
j(I) = (d-1)
!\lim\limits_{q\rightarrow\infty} \frac{\dim_{\k}\Big(\HH_{\mm}^0\left(I^q/I^{q+1}\right)\Big)}{q^{d-1}}.
$$
It was introduced in \cite{ACHILLES_MANARESI_J_MULT} and serves as a generalization of the Hilbert-Samuel multiplicity for non $\mm$-primary ideals.
It is an interesting invariant that has encountered several applications (see ~\cite{FLENNER_O_CARROLL_VOGEL,JMULT_MONOMIAL, JEFFRIES_MONTANO_VARBARO, COMPUTING_J_MULT, POLINI_XIE_J_MULT}).
The next corollary gives a formula for the $j$-multiplicity of any ideal satisfying the conditions of \autoref{thmA}.

\begin{headcor}
	\label{corB}
	Assume all the hypotheses and notations of \autoref{thmA}. 
	Then, the $j$-multiplicity of $I$ is given by
	$$
	j(I) = \frac{1}{2}(n-1)D^{d} \sum_{i=0}^{\lfloor\frac{n-d}{2}\rfloor}\binom{n-2-2i}{d-2}.
	$$
\end{headcor}

In the second main application of \autoref{thmA}, we study rational maps that are determined by a homogeneous minimal set of generators of $I$.
This result adds a new class of rational maps for which the syzygies of the base ideal determine the product of the degrees of the rational map and the corresponding image.
The study of rational maps by using the syzygies of the base ideal is an active research topic (see e.g. ~\cite{AB_INITIO, Simis_cremona,KPU_blowup_fibers,EISENBUD_ULRICH_ROW_IDEALS,Hassanzadeh_Simis_Cremona_Sat,SIMIS_RUSSO_BIRAT,EFFECTIVE_BIGRAD,SIMIS_PAN_JONQUIERES,HASSANZADEH_SIMIS_DEGREES, HULEK_KATZ_SCHREYER_SYZ, MULTPROJ, MULT_SAT_PERF_HT_2, SPECIALIZATION_RAT_MAPS, Kim_Mukundan}).

\begin{headcor}
	\label{corC}
	Assume all the hypotheses and notations of \autoref{thmA}. 
	Choose a homogeneous minimal set of generators $\{f_1,\ldots,f_n\} \subset R$  of $I$ and consider the rational map $\mathcal{F}:\PP_\k^{d-1} \dashrightarrow \PP_\k^{n-1}$ given by
	$$
	\left(x_1:\cdots:x_d\right) \mapsto 	\big(f_1(x_1,\ldots,x_d):\cdots:f_n(x_1,\ldots,x_d)\big).
	$$
	Let $Y \subset \PP_\k^{n-1}$ be the closure of the image of $\mathcal{F}$.
	Then, the following statements hold:
	\begin{enumerate}[(i)]
		\item $\deg(\mathcal{F})\cdot\deg_{\PP_\k^{n-1}}\left(Y\right)=D^{d-1} \sum_{i=0}^{\lfloor\frac{n-d}{2}\rfloor}\binom{n-2-2i}{d-2}.$
		\item $\mathcal{F}$ is birational onto its image if and only if $\deg_{\PP_\k^{n-1}}\left(Y\right)=D^{d-1} \sum_{i=0}^{\lfloor\frac{n-d}{2}\rfloor}\binom{n-2-2i}{d-2}$.
	\end{enumerate}
\end{headcor}

The basic outline of this paper is as follows.
In \autoref{section_mult}, we prove \autoref{thmA}.
In \autoref{section_applications}, we prove \autoref{corB} and \autoref{corC}.

\section{The multiplicity of the saturated special fiber ring}
\label{section_mult}

This section will be divided into two different subsections. 
In the first one, we recall a family of complexes (\cite{KU92}) that will be fundamental in our approach. 
In the second one, we compute the multiplicity of the saturated special fiber ring.

Throughout this paper the following setup will be used.

\begin{setup}
	\label{setup_main}
	Let $\k$ be a field, $R$ be the polynomial ring $R=\k[x_1,\dots,x_d]$ and $I \subset R$ be a homogeneous ideal in $R$.
	Let $\mm \subset R$ be the maximal irrelevant ideal $\mm=\left(x_1,\ldots,x_d\right)$.
	Let $n=\mu(I)$ be the minimal number of generators of $I$.
	Assume the following:
	\begin{enumerate}[(i)]
		\item $I$ is a height $3$ Gorenstein ideal.
		\item  $n \ge d$.
		\item $\varphi$ is a $n \times n$ alternating minimal presentation matrix  of $I$.
		\item Every non-zero entry of the matrix $\varphi$ has degree $D \ge 1$.
		\item $I$ satisfies the condition $G_d$, i.e., $\mu(I_{\pp}) \le \dim(R_{\pp})$  for all $\pp \in V(I) \subset \Spec(R)$  such that $\HT(\pp)<d$.
	\end{enumerate}
\end{setup}

\begin{remark}
	\label{rem_deg_gens}
	Note that the ideal $I$ is minimally generated by $n$ homogeneous polynomials all of the same degree
	$
	\delta := \frac{1}{2}(n-1)D
	$
	(see \cite{BuchsbaumEisenbud}).
\end{remark}

\begin{remark}
	\label{rem_Fitting_conditions}
	In terms of Fitting ideals, $I$ satisfies the condition $G_{d}$ if and only if $\HT(\Fitt_i(I)) > i$ for all $1 \le i < d$.
	So, from the presentation $\varphi$ of $I$, the condition $G_{d}$ is equivalent to $\HT(I_{n-i}(\varphi))>i$ for all $1 \le i < d$. 
	\begin{proof}
		It follows from \cite[Proposition 20.6]{EISEN_COMM}.
	\end{proof}
\end{remark}

The following algebra is the main object of study of this paper, and we shall compute its multiplicity (under the assumptions of \autoref{setup_main}).

\begin{definition}[\cite{MULTPROJ}]
	The saturated special fiber ring of $I$ is given by the graded $\k$-algebra
	$$
	\sF(I) := \bigoplus_{q=0}^\infty {\left[\big(I^q:\mm^\infty\big)\right]}_{q\delta}.
	$$
\end{definition}

We follow an approach very similar to the one used in \cite{MULT_SAT_PERF_HT_2}, and one of our main tools will be to use the complexes $\D_\bullet^q(\varphi)$ (as introduced below in \autoref{subsection_KU_complexes}) in order to obtain approximate resolutions (see, e.g., ~\cite{KPU_DEG_BOUND_LOC_COHOM}, ~\cite{CHARDIN_REG}) of the symmetric powers $\Sym_q\left(I(\delta)\right)$ of $I(\delta)$. 
Another important point in our proof is that, to avoid complicated Hilbert series computations, after some simple reductions we shall use the computation of \cite[Theorem 4.4]{KPU18} (which depends upon the results of \cite{KPU_Hilb_series}).

Let $\AAA$ be the bigraded polynomial ring $\AAA=R[y_1,\ldots,y_n]$ where $\bideg(x_i)=(1,0)$ and $\bideg(y_i)=(0,1)$.
Let $S$ be the standard graded polynomial ring $S=\k[y_1,\ldots,y_n] \subset \AAA$.
The Rees algebra $\Rees(I)=\bigoplus_{q=0}^\infty I^qt^q \subset R[t]$ can be presented as a quotient of $\AAA$ by using the $R$-homomorphism
\begin{eqnarray*}
	\label{presentation_Rees}
	\Psi: \AAA & \longrightarrow & \Rees(I) \subset R[t] \\ \nonumber
	y_i & \mapsto & f_it.
\end{eqnarray*}
We set $\bideg(t)=(-\delta, 1)$, which implies that $\Psi$ is bihomogeneous of degree zero, and so $\Rees(I)$ has a structure of bigraded $\AAA$-algebra.
Note that we immediately obtain the following isomorphism of graded $R$-modules
$$
\Rees(I) \cong \bigoplus_{q=0}^\infty {I^q}(q\delta).
$$
We are mostly interested in the $R$-grading, thus, if $\bbM$ is a bigraded $\AAA$-module and $c \in \ZZ$ a fixed integer, then we write
$$
{\left[\bbM\right]}_c = \bigoplus_{q  \in \ZZ} {\left[\bbM\right]}_{(c,q)}.
$$
Notice that ${\left[\bbM\right]}_c$ has a natural structure as a graded $S$-module.
We recall the following known remark that will allow us to consider certain Hilbert functions.

\begin{remark}
	\label{rem_R_grad_local_cohom}
	For any finitely generated bigraded $\AAA$-module $\bbM$ and any $i \ge 0, c \in \ZZ$, we have that 
	$$
	{\left[\HL^i\left(\bbM\right)\right]}_c 
	$$
	has a natural structure as a finitely generated graded $S$-module (see, e.g., \cite[Theorem 2.1]{CHARDIN_POWERS_IDEALS}, \cite[Proposition 2.7]{MULTPROJ}).
\end{remark}

The special fiber ring of $I$ is defined as $\FF(I):=\Rees(I) / \mm \Rees(I)$.
Since we have 
$$
\Rees(I) = {\left[\Rees(I)\right]}_0\; \bigoplus\;\left( \bigoplus_{c=1}^\infty{\left[\Rees(I)\right]}_c\right)
$$
then we get an isomorphism $\FF(I) \cong {\left[\Rees(I)\right]}_0$ of graded $\k$-algebras.

Following \cite{MULTPROJ}, to study the algebra $\sF(I)$ it is enough to consider the degree zero part in the $R$-grading of the bigraded $\AAA$-module $\HL^1\left(\Rees(I\right))$.

\begin{remark}
	Let $X$ be the scheme $X=\Proj_{R\text{-gr}}\left(\Rees(I)\right)$, where $\Rees(I)$ is only considered as a graded $R$-algebra.
	From \cite[Theorem A4.1]{EISEN_COMM}, we obtain the following short exact sequence
	$$
	0 \rightarrow {\left[\Rees(I)\right]}_0 \rightarrow \HH^0(X, \OO_X) \rightarrow {\left[\HL^1\left(\Rees(I)\right)\right]}_0 \rightarrow 0.
	$$
	By identifying $\FF (I)\cong {\left[\Rees(I)\right]}_0$ and $\sF(I) \cong \HH^0(X, \OO_X)$ (see \cite{MULTPROJ}), we obtain the short exact sequence
	\begin{equation}
	\label{EQ_SPECIAL_AND_SAT}
	0 \rightarrow \FF(I) \rightarrow \sF(I) \rightarrow {\left[\HL^1\left(\Rees(I)\right)\right]}_0 \rightarrow 0.
	\end{equation}
\end{remark}

\begin{remark}
	\label{REM_f.g. Q-mod}
	From \cite[Proposition 2.7$(i)$, Lemma 2.8$(ii)$]{MULTPROJ} we have that $\sF(I)$ and ${\left[\HL^1\left(\Rees(I)\right)\right]}_0$ have natural structures as finitely generated $\mathfrak{F}(I)$-modules.
\end{remark}

As customary, we approximate the Rees algebra with the symmetric algebra by using the following natural short exact sequence
\begin{equation}
\label{EqSymRees}
0 \rightarrow \EEQ \rightarrow \Sym\left(I(\delta)\right) \rightarrow \Rees(I) \rightarrow 0,	
\end{equation}
where $\EEQ$ is the $R$-torsion submodule of $\Sym\left(I(\delta)\right)$.

\begin{lemma}
	\label{LemSymProps}
	The following statements hold:
	\begin{enumerate}[(i)]
		\item $\EEQ=\HL^0\left(\Sym(I(\delta))\right)$.
		\item $\HL^i\left(\Rees(I)\right) \cong \HL^i\left(\Sym(I(\delta))\right)$ for all $i\ge 1$.
	\end{enumerate}
	\begin{proof}
		$(i)$ It follows from \cite[\S 3.7, Eq.~(3.7.3)]{KPU18}.
		
		$(ii)$ Follows identically to \cite[Lemma 2.6$(iii)$]{MULT_SAT_PERF_HT_2}.
	\end{proof}
\end{lemma}

The restriction to the degree zero part in the $R$-grading of the equality $\EEQ=\HL^0\left(\Sym(I(\delta))\right)$ (\autoref{LemSymProps}$(i)$) and the short exact sequence \autoref{EqSymRees} yield the following
\begin{equation}
\label{EQ_SymRees_deg0}
0 \rightarrow {\left[\HL^0\left(\Sym(I(\delta))\right)\right]}_0 \rightarrow S \rightarrow \mathfrak{F}(I) \rightarrow 0,
\end{equation}
under the identifications ${\left[\Sym(I(\delta))\right]}_0 \cong S$ and ${\left[\Rees(I)\right]}_0 \cong \mathfrak{F}(I)$.

\subsection{The complexes of Kustin and Ulrich}
\label{subsection_KU_complexes}
In this short subsection we recall a family of complexes that was introduced in \cite{KU92}.

As in the proof of \cite[Theorem 6.1(b)]{KPU18}, consider the complex
\begin{align*}
\D_\bullet^q(\varphi):\; 0\rightarrow \D_{n-1}^q\rightarrow \D_{n-2}^q\rightarrow\cdots\rightarrow\D_1^q\rightarrow\D_0^q\rightarrow 0
\end{align*}
where
\begin{align}
\label{eq_the_complex}
\D_r^q=\begin{cases}
K_{q-r,r}=R(-rD)^{\beta_r^q} & \text{ if }r\leq\min\{q,n-1\}\\
Q_q=R(-(r-1)D-\frac{1}{2}(n-r+1)D) & \text{ if }r=q+1\leq n-1,q\mathrm{~ odd}\\
0 & \text{ if }r=q+1,q \mathrm{~ even}\\
0 & \text{ if }\min\{ q+2,n\}\leq r.
\end{cases}
\end{align}
The above complex is interesting as the zeroth homology is $\Sym_q\left(I(\delta)\right)$.
Here the symmetric power $\Sym_q\left(I(\delta)\right)$ represents the $q$-th graded piece of $\Sym(I(\delta))$ with respect to the $S$-grading.
 Also, under the assumptions of \autoref{setup_main}, this complex is acyclic on the punctured spectrum of $R$ (\cite[Lemma 6.3]{KPU18}), that is, ${\D_\bullet^q(\varphi)}_\pp$ is acyclic for any $\pp \in \Spec(R) \setminus \{\mm\}$.

To compute the Betti numbers $\beta_r^q$, we use \cite{KU92}. 
Since $\varphi$ is an alternating map, in the definition of \cite[Definition 2.1]{KU92} we have that $\xi=\varphi$ and $\psi=0$, and so $g=\rank(G)=n$ and $f=\rank(F)=0$ (as $F=0$). 
From \cite[Proposition 2.41]{KU92}, we note that 
\begin{align}\label{eq1}
K_{q-r,r}\cong \bigoplus_{v=0}^f L_{q-r+1}^{g-r+v}(G)\otimes_R \wedge^v F
\end{align}
where $L_{q-r+1}^{n-r}(G)$ is the Buchsbaum-Eisenbud free module as defined in \cite[Page 26]{KU92}. 
Hence \autoref{eq1} reduces to 
\begin{align*}
K_{q-r,r}\cong  L_{q-r+1}^{n-r}(G).
\end{align*}
Again, using \cite[Page 26]{KU92} (ranks of $L_p^q$), we have 
\begin{equation}
\label{eq_Betti_numbers}
\beta_r^q={n+q-r\choose n+q-2r}{n+q-2r-1\choose n-r-1}.
\end{equation}

\subsection{The computation of the multiplicity}

During this subsection we compute the multiplicity of the saturated special fiber ring of the ideal $I$.

The following proposition deals with certain local cohomology modules and it will be an important technical tool.

\begin{proposition}
	\label{PropLocCohoSym}
	Let $q \ge 1$.
	Then, we have the following isomorphisms of graded $R$-modules
	$$
	\HH_{i}\Big(\HL^{d}\big(\D_\bullet^q(\varphi)\big)\Big) \cong \begin{cases}
	\HL^{d-i}\left(\Sym_q(I(\delta))\right) \quad \text{ if } i \le d\\
	\HH_{i-d}\left(\D_\bullet^q(\varphi)\right) \qquad\quad\; \text{ if } i \ge d+1,
	\end{cases}	
	$$		
	where $\HL^{d}(\D_\bullet^q(\varphi))$ represents the complex obtained after applying the functor $\HL^{d}(\bullet)$ to $\D_\bullet^q(\varphi)$.
	\begin{proof}
		The main point in the proof is that, from \cite[Lemma 6.3]{KPU18}, the complex $\D_\bullet^q(\varphi)$ is acyclic on the punctured spectrum  of $R$.
		Then, the result follows identically to \cite[Proposition 2.7]{MULT_SAT_PERF_HT_2}.
	\end{proof}
\end{proposition}	

The next lemma contains some dimension computations, and a proof can be found in \cite[Lemma 2.8]{MULT_SAT_PERF_HT_2}.

\begin{lemma}
	\label{LEM_DIM_LOC_SYM_TWO}
	The following statements hold:
	\begin{enumerate}[(i)]
		\item $\dim(\FF(I))=d$.
		\item The corresponding rational $\mathcal{F}:\PP_\k^{d-1} \dashrightarrow \PP_\k^{n-1}$ in \autoref{corC} is generically finite.
		\item $
		\dim\left({\left[\HL^i\left(\Sym(I(\delta))\right)\right]}_0\right) \,\le\, d-1$\; for all $i \ge 2$.
	\end{enumerate}
	\begin{proof}
		See \cite[Lemma 2.8]{MULT_SAT_PERF_HT_2}.
	\end{proof}
\end{lemma}

The following proposition gives a formula which shows that, under the assumptions of \autoref{setup_main}, the multiplicity of $\sF(I)$ only depends on $d$, $n$ and $D$.
This step is fundamental in our approach since it will allow us to use the previous computation of \cite[Theorem 4.4]{KPU18}.

\begin{proposition}
	\label{prop_formular_Hilb_funct}
	Assume \autoref{setup_main}.
	The Hilbert function of $\sF(I)$ is asymptotically a polynomial of degree $d-1$ that satisfies the equation
	$$
	\dim_\k\left({\left[\sF(I)\right]}_q\right) = \binom{q+n-1}{n-1} +{(-1)}^{d-1}\sum_{r=0}^{n-1} {(-1)}^r \binom{rD-1}{d-1}\beta_r^q + O\left(q^{d-2}\right),
	$$
	where $\beta_r^q$ are the Betti numbers of \autoref{eq_Betti_numbers}.
	Therefore, the multiplicity $e\left(\sF(I)\right)$ of $\sF(I)$ only depends on the values of $d$, $n$ and $D$.
	\begin{proof}
		Since $\mathfrak{F}(I) \hookrightarrow \sF(I)$ is an integral extension (see \autoref{REM_f.g. Q-mod}), \autoref{LEM_DIM_LOC_SYM_TWO}$(i)$ implies that
	 	$$
	 	\dim\left(\sF(I)\right)=\dim(\mathfrak{F}(I))=d.
	 	$$
	 	So, it follows that $\dim_\k\left({\left[\sF(I)\right]}_q\right)$ is asymptotically a polynomial of degree $d-1$ (see, e.g., \cite[Theorem 4.1.3]{BRUNS_HERZOG}).
		
		Now, the main idea is to study the homologies of the complexes $\HL^d\left(\D^q_\bullet(\varphi)\right)$.	
		For notational purposes we set 
		$$
		\LL^q_\bullet := \HL^d\left(\D^q_\bullet(\varphi)\right)
		$$	
		for all $q \ge 1$.
		
		Actually, we consider the complexes ${\left[\LL_\bullet^q\right]}_{0}$ of vector spaces over $\k$.
		The additivity of Hilbert functions gives the following equalities
		$$
			\sum_{r=0}^{n-1}{(-1)}^r \dim_\k\Big({\big[\HH_r\left(\LL^q_\bullet\right)\big]}_{0}\Big) = \sum_{r=0}^{n-1} {(-1)}^r \dim_\k\Big({\big[\LL_r^q\big]}_{0}\Big)
		$$
		for $q \ge 1$.
		From \autoref{PropLocCohoSym} and the description \autoref{eq_the_complex} of the complex $\D_\bullet^q(\varphi)$, we obtain that
		$$
			\label{eq_homologies_ge_d+1}
			{\big[\HH_r\left(\LL^q_\bullet\right)\big]}_{0} = {\big[\HH_{r-d}\left(\D_\bullet^q(\varphi)\right)\big]}_0 = 0 \quad \text{ for all } r \ge d+1.
		$$
		Now, for all $r \le d$ we define the function $f_r : \NN \rightarrow \NN$ given by 
		$$
		f_r(q) := \dim_\k\left({\big[\HH_r\left(\LL^q_\bullet\right)\big]}_{0}\right).
		$$
		By using \autoref{PropLocCohoSym} and \autoref{rem_R_grad_local_cohom}, we have that $f_r$ corresponds with the Hilbert function of the finitely generated graded $S$-module
		$
		{\left[\HL^{d-r}\left(\Sym(I(\delta))\right)\right]}_0
		$, that is
		$$
		f_r(q) = \dim_\k\left({\left[\HL^{d-r}\left(\Sym(I(\delta))\right)\right]}_{(0,q)}\right).
		$$ 
		Hence, for all $r \le d$, $f_r(q)$ is a polynomial for $q\gg 0$ and $f_r(q) \in O\left(q^{n-1}\right)$.
		Since 
		$$
		\dim\left({\left[\HL^{d-r}\left(\Sym(I(\delta))\right)\right]}_0\right) \le d-1 \text{ for all } r \le d-2
		$$
		(see \autoref{LEM_DIM_LOC_SYM_TWO}$(iii)$), it follows that $f_r(q) \in O\left(q^{d-2}\right)$ for all $r \le d-2$.
		Therefore, by summing up the above computations we obtain that 
		$$
		{(-1)}^df_d(q) + {(-1)}^{d-1}f_{d-1}(q) = \sum_{r=0}^{n-1} {(-1)}^r \dim_\k\Big({\big[\LL_r^q\big]}_{0}\Big) + O\left(q^{d-2}\right).
		$$
		Note that for $q\gg 0$ we have that $\D_r^q=R(-rD)^{\beta_r^q}$ (see \autoref{eq_the_complex}).
		Thus, the isomorphism 
		$$
		\HL^d(R) \cong \frac{1}{x_1\cdots x_d}\k[x_1^{-1},\ldots,x_d^{-1}]
		$$
		yields that $\dim_\k\Big({\big[\LL_r^q\big]}_{0}\Big)=\binom{rD-1}{d-1}\beta_r^q$.
		So, for $q\gg 0$, we get that 
		\begin{equation}
			\label{eq_euler_f_i_s}
			{(-1)}^df_d(q) + {(-1)}^{d-1}f_{d-1}(q) = \sum_{r=0}^{n-1} {(-1)}^r \binom{rD-1}{d-1}\beta_r^q + O\left(q^{d-2}\right).
		\end{equation}
		
		The short exact sequence in \autoref{EQ_SPECIAL_AND_SAT} and \autoref{LemSymProps}$(ii)$ give 
		$$
		\dim_\k\left({\left[\sF(I)\right]}_q\right) = \dim_\k\left({\left[\HL^1\left(\Sym(I(\delta))\right)\right]}_{(0,q)}\right) + \dim_\k\left({\left[\mathfrak{F}(I)\right]}_q\right),
		$$
		and then from the short exact sequence in \autoref{EQ_SymRees_deg0} we have that 
		\begin{align}
			\label{eq_Hilb_funct_sF_in_f_i_s}
			\begin{split}
			\dim_\k\left({\left[\sF(I)\right]}_q\right) &= \dim_\k\left({\left[\HL^1\left(\Sym(I(\delta))\right)\right]}_{(0,q)}\right) + \dim_\k\left({\left[S\right]}_q\right) \\
			& \qquad\qquad\qquad\qquad\qquad- \dim_\k\left({\left[\HL^0\left(\Sym(I(\delta))\right)\right]}_{(0,q)}\right)\\
			&= f_{d-1}(q) + \binom{q+n-1}{n-1} - f_d(q).
			\end{split}
		\end{align}
		Finally, by combining \autoref{eq_euler_f_i_s} and \autoref{eq_Hilb_funct_sF_in_f_i_s}, for $q \gg 0$ we obtain the equation
		$$
		\dim_\k\left({\left[\sF(I)\right]}_q\right) = \binom{q+n-1}{n-1} +{(-1)}^{d-1}\sum_{r=0}^{n-1} {(-1)}^r \binom{rD-1}{d-1}\beta_r^q + O\left(q^{d-2}\right).
		$$
		So, the result follows.
	\end{proof}		
\end{proposition}

The lemma below provides lower bounds for the grade of certain determinantal ideals. 
It will be used to obtain ideals that satisfy the conditions of \autoref{setup_main} and that are convenient to apply the computation of \cite[Theorem 4.4]{KPU18}.
It is likely part of the folklore, but we include a proof for the sake of completeness; for a similar setup, see \cite[Proposition 4.2]{SPECIALIZATION_RAT_MAPS}.

\begin{lemma}
	\label{lem_det_grade_lb}
	Let $m \ge 3$ be an odd integer.
	Let $\mathbf{z}=(z_{i,j,k})$ be a new set of variables with $1\le i \le m$, $i+1 \le j \le m$ and $1 \le k \le d$.
	Let $T$ be the polynomial ring $T=R[\mathbf{z}]$.
	Let $M$ be the alternating $m\times m$ matrix with entries in $T$ given by 
	$$
	M = \left( \begin{array}{ccccc}
	0 & p_{1,2} & \cdots & p_{1,m-1} & p_{1,m} \\
	-p_{1,2} & 0 & \cdots & p_{2,m-1 }& p_{2,m}\\
	\vdots & \vdots & & \vdots &\vdots\\
	-p_{1,m-1} & -p_{2,m-1} & \cdots & 0 & p_{m-1,m}\\
	-p_{1,m} & -p_{2,m} & \cdots & -p_{m-1,m} & 0\\
	\end{array}
	\right)	
	$$	
	where each polynomial $p_{i,j} \in T$ is given by 
	$$
	p_{i,j} = x_1^Dz_{i,j,1}  + x_2^Dz_{i,j,2}  + \cdots + x_d^Dz_{i,j,d}.
	$$
	Then 
	$$
	\grade\big(I_t(M)\big) \ge 
	\begin{cases}
		\min\{ m+2-t, d \}  \qquad \text{if $t$ is even} \\
		\min\{ m+1-t, d \}  \qquad \text{if $t$ is odd}
	\end{cases}
	$$
	for $1 \le t \le m-1$.
	\begin{proof}
		For any ideal $J \subset T$, denote by $\text{Rad}\left(J\right)$ its radical ideal.
		Since $\text{Rad}\left(I_{t-1}(M)\right) = \text{Rad}\left(I_t(M)\right)$ when $t$ is even (see, e.g., \cite[Corollary 2.6]{BuchsbaumEisenbud}), it suffices to show that $\grade(I_t(M)) \ge \min\{ m+2-t, d \}$ when $t$ is even.
		
		We proceed by induction on $t$ and assuming that $t$ is even with $2 \le t \le m-1$.
		The case $t=2$ follows from \cite[Lemma 4.1]{SPECIALIZATION_RAT_MAPS} since $\text{Rad}(I_2(M))=\text{Rad}(I_1(M))$ and $I_1(M)$ is generated by the $p_{i,j}$'s themselves.
		
		Now suppose that $2 < t \le m-1$.
		Let $Q \in \Spec(T)$ be a prime ideal containing $I_t(M)$.
		If $Q$ contains all the polynomials $p_{i,j}$, then again \cite[Lemma 4.1]{SPECIALIZATION_RAT_MAPS} yields 
		$$
		\depth(T_{Q}) \ge \min\Big\{\frac{m(m-1)}{2},d\Big\} \ge \min\{m+2-t,d\}.
		$$
		Otherwise, say, $p_{m-1,m} \not\in Q$.
		
		Let $M^\prime$ denote the $(m-2)\times(m-2)$ submatrix of $M$ of the first $m-2$ rows and first $m-2$ columns. 
		Clearly, 
		$$
		I_{t-2}\left(M^\prime\right)T_{p_{m-1,m}} \,\subset\, I_{t}\left(M\right)T_{p_{m-1,m}}
		$$
		in the localization  $T_{p_{m-1,m}}$.
		The inductive hypothesis gives 
		\begin{align*}
		\depth(T_Q) &\ge 
		\grade\left(I_{t}\left(M\right)T_{p_{m-1,m}}\right)	\\
		&\ge \grade\left(I_{t-2}\left(M^\prime\right)T_{p_{m-1,m}}\right)\\
		&\ge 	\grade\left(I_{t-2}\left(M^\prime\right)\right) \\
		&\ge \min\{(m-2)+2-(t-2),d\}.
		\end{align*}
		Therefore, $\depth(T_Q)\ge \min\{m+2-t,d\}$, and so the result follows.
	\end{proof}
\end{lemma}

Now, we are ready for the main result of this paper. 

\begin{proof}[Proof of \autoref{thmA}]
	From \autoref{prop_formular_Hilb_funct}, $e\left(\sF(I)\right)$ depends only on $d$, $n$ and $D$, and so we may substitute $I$ by a suitable ideal with the same numerics.
	
	Without any loss of generality, assume that $\k$ is an infinite field.
	Let $M \in {R[\mathbf{z}]}^{n \times n}$ be the matrix of \autoref{lem_det_grade_lb} with $m=n$.
	Then, by using \cite[Lemma 3.4]{MOREY_ULRICH} and \cite[Lemma 2.1]{HUNEKE_ULRICH}, there exists a dense open subset $U \subset \k^{d\frac{n(n-1)}{2}}$ such that, for any $\underline{\alpha} = (\alpha_{i,j,k}) \in U$, the specialization $\Pi_{\underline{\alpha}}(M) \in R^{n\times n}$ obtained by setting $z_{i,j,k} \mapsto \alpha_{i,j,k}$ satisfies the following inequalities 
	\begin{equation}
		\label{eq_ineq_grade_sepcializations}
		\grade\big(I_{t}\left(\Pi_{\underline{\alpha}}(M)\right)\big) \ge \grade\left(I_t\left(M\right)\right) \ge 
		\begin{cases}
		\min\{ n+2-t, d \}  \qquad \text{if $t$ is even} \\
		\min\{ n+1-t, d \}  \qquad \text{if $t$ is odd}
		\end{cases}
	\end{equation}
	for $1 \le t \le n-1$.

	Fix any $\underline{\alpha} = (\alpha_{i,j,k}) \in \k^{d\frac{n(n-1)}{2}}$ for which the inequalities of \autoref{eq_ineq_grade_sepcializations} hold, and set $\varphi = \Pi_{\underline{\alpha}}(M) \in R^{n\times n}$ and $I=\text{Pf}_{n-1}(\varphi) \subset R$.	
	Hence, \cite[Theorem 2.1, Corollary 2.6]{BuchsbaumEisenbud}, \autoref{rem_Fitting_conditions} and \autoref{eq_ineq_grade_sepcializations} imply that $I$ satisfies all the conditions of \autoref{setup_main}.
	
	Note that all the entries of $\varphi$ are $\k$-linear combinations of the monomials $x_1^D, \ldots, x_d^D$.
	Let $R^\prime$ be the polynomial ring $R^\prime=\k[w_1,\ldots,w_d]$ and $\Phi:R^\prime \rightarrow R$ be the $\k$-algebra homomorphism given by $\Phi(w_i) = x_i^D$.
	Set $\varphi^\prime=\Phi^{-1}\left(\varphi\right) \in {R^\prime}^{n\times n}$ and $I^\prime = \text{Pf}_{n-1}\left(\varphi^\prime\right) \in R^\prime$.
	Thus, $I^\prime$ is a linearly presented Gorenstein ideal of height $3$ that satisfies the condition $G_d$.
	Let $f_1^\prime, \ldots, f_n^\prime \in R^\prime$ and $f_1, \ldots, f_n \in R$ be the $(n-1)\times(n-1)$ pfaffians of $\varphi^\prime$ and $\varphi$, respectively; note that $f_i=\Phi(f_i^\prime)$.
	
	Let $\mathcal{F} : \PP_\k^{d-1}\cong\Proj(R) \dashrightarrow  \PP_\k^{n-1}$, $\mathcal{F}^\prime : \PP_\k^{d-1}\cong\Proj\left(R^\prime\right) \dashrightarrow  \PP_\k^{n-1}$ and $\mathcal{G} : \PP_\k^{d-1}\cong\Proj\left(R\right) \rightarrow \PP_\k^{d-1}\cong\Proj\left(R^\prime\right)$ be rational maps with representatives $\left(f_1:\cdots:f_n\right)$, $\left(f_1^\prime:\cdots:f_n^\prime\right)$ and $\left(x_1^D:\cdots:x_d^D\right)$, respectively.
	Then, we obtain the following commutative diagram
		\begin{center}		
		\begin{tikzpicture}
		\matrix (m) [matrix of math nodes,row sep=2.6em,column sep=6em,minimum width=2em, text height=1.5ex, text depth=0.25ex]
		{
			\PP_\k^{d-1} \cong \Proj\left(R\right) & & \PP_\k^{n-1}\\
			 & \PP_\k^{d-1} \cong \Proj\left(R^\prime\right)  &  \\
		};
		\path[-stealth]
		(m-1-1) edge node [above] {$\mathcal{G}$} (m-2-2)
		;				
		\draw[->,dashed] (m-1-1)--(m-1-3) node [midway,above] {$\mathcal{F}$};	
		\draw[->,dashed] (m-2-2)--(m-1-3) node [midway,above] {$\mathcal{F}^\prime$};	
		\end{tikzpicture}			
	\end{center}
	and so it follows that $\deg\left(\mathcal{F}\right)=\deg\left(\mathcal{F}^\prime\right)\cdot\deg(\mathcal{G})$.
	Let $Y \subset \PP_\k^{n-1}$ be the closure of the image of $\mathcal{F}$.
	Note that $Y$ coincides with the closure of the image of $\mathcal{F}^\prime$.
	From \cite[Proposition 4.1, Theorem 4.4]{KPU18} we get that $\deg(\mathcal{F}^\prime)=1$ and 
	$
	\deg_{\PP_\k^{n-1}}(Y)=\sum_{i=0}^{\lfloor\frac{n-d}{2}\rfloor}\binom{n-2-2i}{d-2}.
	$
	On the other hand, it is well-known that $\deg(\mathcal{G})=D^{d-1}$ (see, e.g., \cite[Observation 3.2]{KPU_blowup_fibers}, \cite[Theorem 3.3]{MULTPROJ}).
	
	Finally, \cite[Theorem 2.4]{MULTPROJ} yields the equality
	$$
	e\left(\sF(I)\right)=\deg(\mathcal{F})\cdot\deg_{\PP_\k^{n-1}}\left(Y\right)=D^{d-1}\sum_{i=0}^{\lfloor\frac{n-d}{2}\rfloor}\binom{n-2-2i}{d-2}.
	$$
	So, the result follows.
\end{proof}

\section{Some applications}
\label{section_applications}

Throughout this section, we continue using \autoref{setup_main}.

Here we provide some applications that follow straightforwardly from the computation of \autoref{thmA} and \cite{MULTPROJ}.
More specifically, under the assumptions of \autoref{setup_main}, we give an exact formula for the $j$-multiplicity of $I$ and we study a rational map 
\begin{align}
\label{EQ_RAT_MAP}
&\mathcal{F} : \PP_\k^{d-1}\dashrightarrow \PP_\k^{n-1}\\
\nonumber
\left(x_1:\cdots:x_d\right) &\mapsto 	\big(f_1(x_1,\ldots,x_d):\cdots:f_n(x_1,\ldots,x_d)\big)
\end{align}
that corresponds to a homogeneous minimal set of generators $\{f_1,\ldots,f_n\} \subset R$ of $I$.
Let $Y \subset \PP_\k^{n-1}$ be the closure of the image of $\mathcal{F}$, and $K(\PP_\k^{d-1})$ and $K(Y)$ be the fields of rational functions of $\PP_\k^{d-1}$ and $Y$, respectively.
The degree of $\mathcal{F}$ is equal to the degree of the field extension 
$$
\deg(\mathcal{F}) = \left[K(\PP_\k^{d-1}):K(Y)\right].
$$

First, we have a formula for the $j$-multiplicity of $I$.

\begin{proof}[Proof of \autoref{corB}]
	From \cite[Lemma 2.10]{MULTPROJ} we know that $j(I)=\delta\cdot e\left(\sF(I)\right)$.
	Then, \autoref{rem_deg_gens} and \autoref{thmA} yield
	$$
	j(I)=\frac{1}{2}(n-1)D\cdot D^{d-1} \sum_{i=0}^{\lfloor\frac{n-d}{2}\rfloor}\binom{n-2-2i}{d-2}=\frac{1}{2}(n-1) D^{d} \sum_{i=0}^{\lfloor\frac{n-d}{2}\rfloor}\binom{n-2-2i}{d-2},
	$$
	and so the result follows.
\end{proof}

Now, by using the formula of \autoref{thmA}, we obtain some important consequences for the rational map $\mathcal{F}$ in \autoref{EQ_RAT_MAP}.

\begin{proof}[Proof of \autoref{corC}]
	From \cite[Theorem 2.4]{MULTPROJ} we have that 
	$$
	e\left(\sF(I)\right)=\deg(\mathcal{F})\cdot \deg_{\PP_\k^{n-1}}(Y).
	$$
	So, the result is clear from \autoref{thmA}.
\end{proof}

\section*{Acknowledgments}
We thank the reviewer for his/her suggestions for the improvement of this work.
The use of \textit{Macaulay2} \cite{MACAULAY2} was very important in the preparation of this paper.

\bibliographystyle{elsarticle-num} 
%\bibliography{references}

\end{document}